\documentclass{ifacconf}
\usepackage{graphicx}
\usepackage{natbib}
\usepackage{xcolor}
\usepackage{amsfonts}
\usepackage{amsmath}
\usepackage{algorithm}
\usepackage[noend]{algpseudocode}
\usepackage{url}

\begin{document}
\begin{frontmatter}

\title{An adaptive model hierarchy for data-augmented training of kernel models for reactive flow\thanksref{footnoteinfo}}

\thanks[footnoteinfo]{%
  Funded by BMBF under contracts 05M20PMA and 05M20VSA.
  Funded by the Deutsche Forschungsgemeinschaft (DFG, German Research Foundation) under contracts OH 98/11-1 and SCHI 1493/1-1, as well as under Germany's Excellence Strategy EXC 2044 390685587, Mathematics Münster: Dynamics -- Geometry -- Structure, and EXC 2075 390740016, Stuttgart Center for Simulation Science (SimTech).
}

\author[First]{B.~Haasdonk}
\author[Second]{M.~Ohlberger}
\author[Second]{F.~Schindler}

\address[First]{Institute of Applied Analysis and Numerical Simulation, Pfaffenwaldring 57, D-70569 Stuttgart (e-mail: \textnormal{\texttt{haasdonk@mathematik.uni-stuttgart.de}}).}
\address[Second]{Mathematics Münster, Westfälische Wilhelms-Universität Münster, Einsteinstr.~62, D-48149 Münster (e-mail: \textnormal{\texttt{\{mario.ohlberger,felix.schindler\}@uni-muenster.de}})}

\end{frontmatter}

\section{Reference model}
\label{sec:reference_model}

We are interested in constructing efficient and accurate models to approximate time-dependent quantities of interest (QoI) $f \in L^2\big(\mathcal{P}; L^2([0, T])\big) $ in the context of reactive flow, with $T > 0$ and where $\mathcal{P} \subset \mathbb{R}^p$ for $p \geq 1$ denotes the set of possible input parameters.
As a class of QoI functions, we consider those obtained by applying linear functionals $s_\mu \in V'$ to solution trajectories $c_\mu \in L^2(0, T; V)$ of, e.g., parametric parabolic partial differential equations.
Thus, $f(\mu; t) := s_\mu(c_\mu(t))$, where
for each parameter $\mu \in \mathcal{P}$, the concentration $c_\mu$ with $\partial_t c_\mu \in L^2(0, T; V')$ and initial condition $c_0 \in V$ is the unique weak solution of
\begin{align}
  \left<\partial_t c_\mu, v\right> + a_\mu(c_\mu, v) = l_\mu(v) &&\forall\; v \in V,\;c_\mu(0) = c_0.
  \label{eq:pde}
\end{align}
Here, $V \subset H^1(\Omega) \subset L^2(\Omega) \subset V'$ denotes a Gelfand triple of Hilbert-spaces associated with a spatial Lipschitz-domain $\Omega$ and, for $\mu \in \mathcal{P}$, $l_\mu \in V'$ denotes a continuous linear functional and $a_\mu: V \times V \to \mathbb{R}$ a continuous coercive bilinear form.

As a basic model for reactive flow in catalytic filters, \eqref{eq:pde} could stem from a single-phase one-dimensional linear advection-diffusion-reaction problem with Dammköhler- and P\'{e}clet-numbers as input (thus $p = 2$), where $c$ models the dimensionless molar concentration of a species and the break-through curve $s$ measures the concentration at the outflow, as detailed in \cite{GHI+2021}.

Since direct evaluations of $f$ are not available, we resort to a full order model (FOM) as reference model, yielding
\begin{align}
  f_h: \mathcal{P} \to \mathbb{R}^{N_T}\text{ for }N_T \geq 1, &&f_h(\mu; t) := s_\mu(c_{h, \mu}(t)),
  \label{eq:f_h}
\end{align}
which we assume to be a sufficiently accurate approximation of the QoI.
For simplicity, we consider a $P^1$-conforming Finite Element space $V_h \subset V$ and obtain the FOM solution trajectory $c_{h, \mu} \in L^2(0, T; V_h)$ by Galerkin projection of \eqref{eq:pde} onto $V_h$ and an implicit Euler approximation of the temporal derivative.

\section{Surrogate models}

The evaluation of \eqref{eq:f_h} may be arbitrarily costly, in particular in multi- or large-scale scenarios where $\dim V_h \gg 1$, but also if $N_T \gg 1$ due to long-time integration or when a high resolution of $f_h$ is required.
We thus seek to build a machine learning (ML) based surrogate model
\begin{align}
  f_\text{ml}: \mathcal{P} \to \mathbb{R}^{N_T},\;f_\text{ml}(\mu; t_n) \approx f_h(\mu; t_n)\;\forall 1\leq n \leq N_T,
  \label{eq:f_ml}
\end{align}
to predict all values $\{f_\text{ml}(\mu; t_n)\}_{n=1}^{N_T}$ at once, without time-integration.
Such models based on Neural Networks or Kernels typically rely on a large amount of training data
\begin{align}
  \big\{\big(\mu, f_h(\mu)\big)\,\big|\,\mu \in \mathcal{P}_\text{train}\big\},&& \mathcal{P}_\text{train}\subset\mathcal{P},&& |\mathcal{P}_\text{train}|\gg 1,
\end{align}
rendering their training prohibitively expensive in the aforementioned scenarios; we refer to \cite{GHI+2021} and the references therein and in particular to \cite{SH2019}.
In \cite{GHI+2021} we thus seek to employ an intermediate surrogate to generate sufficient training data.

\subsection{Structure preserving Reduced Basis models}

The idea of projection-based model order reduction by Reduced Basis (RB) methods is to approximate the state $c_h$ in a low-dimensional subspace $V_\text{rb} \subset V_h$ and to obtain online-efficient approximations of $f_h$ by Galerkin projection of the FOM detailed in Section \ref{sec:reference_model} onto $V_\text{rb}$ and a pre-computation of all quantities involving $V_h$ in a possibly expensive offline-computation; we refer to \cite{MRS2016} and the references therein.
Using such structure preserving reduced order models (ROM)s we obtain RB trajectories $c_{\text{rb},\mu} \in L^2(0, T; V_\text{rb})$ and a RB model
\begin{align}
  f_\text{rb}: \mathcal{P} \to \mathbb{R}^{N_T}, &&f_\text{rb}(\mu; t) := s_\mu(c_{\text{rb}, \mu}(t)),
  \label{eq:f_rb}
\end{align}
with a computational complexity independent of $\dim V_h$, the solution of which, however, still requires time-in\-te\-gra\-tion.

The quality and efficiency of RB models hinges on the problem adapted RB space $V_\text{rb}$ which could be constructed in an iterative manner steered by a posteriori error estimates using the POD-greedy algorithm from \cite{Haa2013}.
Instead, we obtain by the method of snapshots
\begin{align}
  V_\text{rb} := \left< \texttt{POD(}\{c_{h, \mu} \,|\, \mu \in \mathcal{P}_\text{rb} \texttt{)} \right>,&& \text{with }\mathcal{P}_\text{rb} \subset \mathcal{P}
\end{align}
consisting of only few a priori selected parameters (e.g.~the outermost four points in $\mathcal{P}$), where we use the hierarchic approximate POD from \cite{HLR2018} for $N_T \gg 1$ to avoid computing the SVD of a dense snapshot Gramian of size $N_T^2$.

\subsection{Kernel models}

While still requiring time-integration, we can afford to use RB ROMs to generate a sufficient amount of training data
\begin{align*}
  X_\text{train} = \big\{\big(\mu, f_\text{rb}(\mu)\big)\,\big|\,\mu \in \mathcal{P}_\text{ml}\big\} \cup \big\{\big(\mu, f_h(\mu)\,\big|\,\mu\in\mathcal{P}_\text{rb}\big)\big\},
\end{align*}
augmented by the FOM-data available as a side-effect from generating $V_\text{rb}$.
Using this data, we obtain the ML model $f_\text{ml}$ from \eqref{eq:f_ml} using the vectorial greedy orthogonal kernel algorithm from \cite{SH2019}.

While resulting in substantial computational gains, the presented approach from \cite{GHI+2021} still relies on the traditional offline/online splitting of the computational process to train the RB ROM as well as the ML model to be valid for all of $\mathcal{P}$, requiring a priori choices regarding $\mathcal{P}_\text{rb}$ and $\mathcal{P}_\text{ml}$ with a significant impact on the overall performance and applicability of these models.

\section{An adaptive model hierarchy}

\cite{KMOSV2021} introduced an approach beyond the classical offline/online splitting where a RB ROM is adaptively enriched based on rigorous a posteriori error estimates, following the path of an optimization procedure through the parameter space.
Similarly, we propose an adaptive enrichment yielding a hierarchy of FOM, RB ROM and ML models, based on the standard residual-based a posteriori estimate on the RB output error, $\|f_h(\mu) - f_\text{rb}(\mu)\|_{L^2([0, T]))} \leq \Delta_\text{rb}(\mu)$, for which we refer to the references in \cite{MRS2016}.

\begin{algorithm}
  \caption{Adaptive QoI model generation}
  \label{alg:adaptive_model}
  \begin{algorithmic}[1]
    \Require ROM tolerance $\varepsilon>0$, ML trust/train criteria
    \State $X_\text{train} = \emptyset$, $\Phi_\text{RB} = \{\}$, $V_\text{rb} := \left<\Phi_\text{rb}\right>$, $f_\text{ml} := 0$
    \ForAll{$\mu \in \mathcal{P}$ selected by outer loop}
      \If{ML model is trustworthy}
        \Return $f_\text{ml}(\mu)$
      \Else
        \State compute $f_\text{rb}(\mu)$, $\Delta_\text{rb}(\mu)$
        \If{$\Delta_\text{rb}(\mu) \leq \varepsilon$}
          \State collect $X_\text{train} = X_\text{train} \cup \{(\mu, f_\text{rb}(\mu))\}$
          \State (optionally) fit ML model,
          \Return $f_\text{rb}(\mu)$
          \Else{\hfill{{\color{gray}{$\Pi_{\Phi_\text{rb}}$: orth.~proj. onto $\left<\Phi_\text{rb}\right>$}}}}
          \State compute $f_h(\mu)$
          \State enrich $\Phi_\text{rb} = \Phi_\text{rb}\cup\texttt{POD(}c_{h}(\mu)-\Pi_{\Phi_\text{rb}}[c_h(\mu)]\texttt{)}$
          \State update RB ROM
          \State collect $X_\text{train} = X_\text{train} \cup \{(\mu, f_h(\mu))\}$
          \State (optionally) fit ML model,
          \Return $f_h(\mu)$
        \EndIf
      \EndIf
    \EndFor
  \end{algorithmic}
\end{algorithm}

As a means to judge if a ML model is trustworthy, we propose a manual validation using the following a posteriori error estimate on the ML QoI error.
While not as cheaply computable as $f_\text{ml}$, it still allows to validate the ML model without computing $f_h$.

\begin{prop}[ML model a posteriori error estimate]
  \mbox{}\\Let $f_\text{rb}(\mu)$, $f_\text{ml}(\mu) \in \mathbb{R}^{N_T}$ denote the RB ROM and ML model approximations of $f_h(\mu)$, respectively, and
  let $\Delta_\text{rb}(\mu)$ denote an upper bound on the RB-output error.
  We then have by triangle inequality for all $\mu \in \mathcal{P}$
  \begin{align}
    \|f_h(\mu) - f_\text{ml}(\mu)&\|_{L^2([0, T])} \leq \Delta_\text{rb}(\mu)\\\notag
               &+ \|f_\text{rb}(\mu) - f_\text{ml}(\mu)\|_{L^2([0, T])},
  \end{align}
  where the right hand side is computable with a computational complexity independent of $\dim V_h$.
\end{prop}

\begin{figure}[h]
\begin{center}
\includegraphics[width=8.4cm,height=3.6cm]{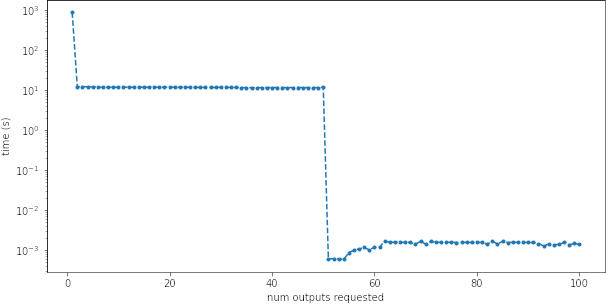}
\caption{Each dot correspond to the input-to-output query time of the adaptive model from Algorithm \ref{alg:adaptive_model} applied to \cite{GHI+2021}.
}
\label{fig:results}
\end{center}
\end{figure}

Applying Algorithm \ref{alg:adaptive_model} to the example of one-dimensional single-phase reactive flow from the last row of Table 1 in \cite{GHI+2021}, with $\dim \mathcal{P}=2$, $N_T = 24576$ time steps, $\dim V_h = 65537$ gives the behaviour shown in Figure \ref{fig:results}, where we set $\varepsilon = 10^{-2}$, retrain the ML model every 10 collected samples and unconditionally trust the ML model as soon as $|X_\text{train}| \geq 50$.\footnote{%
  The experiments were performed using \texttt{pyMOR} from \cite{MRS2016} and \texttt{dune-gdt} from \url{https://docs.dune-gdt.org/}.
}
For the considered diffusion dominated regime, we only require a single evaluation of $f_h$ (yielding a $\dim V_\text{rb} = 15$-dimensional RB ROM), which results in even further computational savings, compared to the results obtained in \cite{GHI+2021}.

\end{document}